\documentclass[12pt]{article}
\usepackage{amsmath,amssymb,amscd,amsthm,latexsym,amsfonts}

\begin{document}

\newtheorem{question}{Question}
\newtheorem{conjecture}[equation]{Conjecture} 
\newtheorem{theorem}[equation]{Theorem}
\newtheorem{lemma}[equation]{Lemma}
\newtheorem*{claim}{Claim}
\newtheorem{corollary}[equation]{Corollary}
\newtheorem{prop}[equation]{Proposition}
\newtheorem{problem}[equation]{Problem}

\theoremstyle{remark}
\newtheorem*{remark}{Remark}

\newcommand{\dv}{\operatorname{div}}
\newcommand{\R}{\operatorname{Re}}
\newcommand{\supp}{\operatorname{supp}}
\newcommand{\dist}{\operatorname{dist}}
\newcommand{\Lip}{\operatorname{Lip}}
\newcommand{\diam}{\operatorname{diam}}
\newcommand{\epi}{\operatorname{Epi}}

\newcommand{\NN}{\mathbb{N}}
\newcommand{\RR}{\mathbb{R}}
\newcommand{\ZZ}{\mathbb{Z}}
\newcommand{\QQ}{\mathbb{Q}}
\newcommand{\CC}{\mathbb{C}}
\newcommand{\KK}{\mathbb{K}}
\newcommand{\rn}{\RR^n}
\newcommand{\srn}{{\scriptscriptstyle \RR}^n}
\newcommand{\rplus}{\RR_+}
\newcommand{\rplusbar}{\ol{\RR_+}}

\newcommand{\ep}{\varepsilon}
\newcommand{\si}{\sigma}
\newcommand{\dxdtovert}{\frac{dxdt}{t}}
\newcommand{\comp}{{}^{\textnormal{c}}}

\newcommand{\calA}{\mathcal{A}}
\newcommand{\calB}{\mathcal{B}}
\newcommand{\calC}{\mathcal{C}}
\newcommand{\calD}{\mathcal{D}}
\newcommand{\calE}{\mathcal{E}}
\newcommand{\calF}{\mathcal{F}}
\newcommand{\calG}{\mathcal{G}} 
\newcommand{\calH}{\mathcal{H}} 
\newcommand{\calI}{\mathcal{I}} 
\newcommand{\calJ}{\mathcal{J}} 
\newcommand{\calK}{\mathcal{K}} 
\newcommand{\calL}{\mathcal{L}}
\newcommand{\calM}{\mathcal{M}}
\newcommand{\calN}{\mathcal{N}}
\newcommand{\calO}{\mathcal{O}}
\newcommand{\calP}{\mathcal{P}}
\newcommand{\calQ}{\mathcal{Q}}
\newcommand{\calR}{\mathcal{R}}
\newcommand{\calS}{\mathcal{S}}
\newcommand{\calT}{\mathcal{T}}
\newcommand{\calU}{\mathcal{U}}
\newcommand{\calV}{\mathcal{V}}
\newcommand{\calW}{\mathcal{W}}
\newcommand{\calX}{\mathcal{X}}
\newcommand{\calY}{\mathcal{Y}}
\newcommand{\calZ}{\mathcal{Z}}

\title{Lectures on the Kato square root problem}

\author{Pascal Auscher\\LAMFA, CNRS FRE 2270\\Universit\'e de
Picardie-Jules Verne\\33, rue Saint Leu\\F-80039
Amiens Cedex 1\\email:
auscher@mathinfo.u-picardie.fr
% \and Steve
%Hofmann\\Department of Mathematics\\University of
%Missouri-Columbia\\Columbia, MO 65211\\email:
%hofmann@math.missouri.edu
} 

\date{July 31, 2001}
\maketitle

\begin{abstract} This is the text of a series of three
lectures  given at the CMA of the Australian National
University on the recent solution of the  square root
problem for divergence form elliptic operators, a
long-standing conjecture posed by Kato in the early 60's. In
this text,  the motivations for this problem and its
situation are given.  The ideas from harmonic analysis on the
T(1) theorem and T(b) theorem for square functions are
described. In particular,  an apparently new
formulation of a local T(b) theorem for square functions is
stated. The ideas of the full proof are presented.   I want to
thank the CMA at the Australian National University for
inviting me during the special program on scattering theory
and spectral problems   and for the nice and stimulating
atmosphere  created by the  mathematicians at the CMA.  
\end{abstract}
\vfill\break
\tableofcontents
\vfill\break  

\section{Elliptic operators}

Consider an open subset $\Omega$ of $\RR^n$, $n\ge 1$. Let
$V$ be a closed subspace of $H^m(\Omega)=W^{m,2}(\Omega)$
which contains $H^m_0(\Omega)$, the closure of smooth
functions supported in $\Omega$ in $H^m(\Omega)$.

Let $N,m$ be positive integers and define a sesquilinear form
on $V^N\times V^N$ by
$$
Q(f,g) = \int_\Omega  \sum_{|\alpha|,|\beta| \le m
\atop 1 \le i, j \le N} 
 a_{\alpha\beta}^{ij}(x) \partial^\beta f_j(x) 
\partial^\alpha\overline g_i(x) \, dx
$$
Here ${f}= (f_1, \ldots ,
f_N)$  and ${g}= (g_1, \ldots ,
g_N)$ belongs to $V^N$, and the  coefficients 
$a_{\alpha\beta}^{ij}$ are   complex-valued  $L^\infty$ 
functions on $\Omega$. We  use the standard notations of
 differential calculus in $\RR^n$: multiindices, partials...

One assumes that 
 \begin{equation}\label{eqLhombounded}
| Q(f,g)| \le \Lambda \| \nabla^m f\|_2 \|
\nabla^m g\|_2 + \kappa'\|f\|_2\|g\|_2
\end{equation}
and the G\aa rding inequality   
\begin{equation}\label{eqgarding}
\R Q(f,f) \ge \lambda \| \nabla^m f\|_2^2 -\kappa\|f\|_2^2
\end{equation}
for some  $\lambda>0$, $\kappa,\kappa'\ge 0$ and $\Lambda
<+\infty$ independent of 
$f,g
\in  V^N$. Here, $\nabla^kf= (\partial^\alpha
f)_{|\alpha|=k}$  and $\| \nabla^k
f\|_2=
\big(\sum_{|\alpha|= k}
\int_{\Omega}|\partial^\alpha f|^2
\big)^{1/2}$....  

A well-known representation theorem of Kato asserts that one
can represent the form by 
$$
Q(f,g) = \langle Lf, g\rangle, \quad f\in \calD(L), g\in V^N
$$
where $D(L)$ is the subspace of those $f\in V^N$ such that 
$g\mapsto Q(f,g)$, originally defined on $V^N$, extends to a
bounded anti-linear form on $L^2(\Omega, \CC^N)$. As usual,
it is convenient to denote the operator (system) as
\begin{eqnarray}\label{eqsystem}
(L{f})_i= \sum_{|\alpha|,|\beta| \le m \atop 1 \le j \le N}
(-1)^{|\alpha|}
\partial^\alpha (a_{\alpha\beta}^{ij} \partial^\beta f_j ), \quad
1\le i\le N.
\end{eqnarray}
In fact,  the operator $L$ is defined from
$V^N$ into its dual and $\calD(L)$  can be seen as the
subspace of
$f\in V^N$ such that $Lf \in L^2(\Omega,\CC^N)$.
The restriction of
$L$ to
$\calD(L)$ is a maximal-accretive operator and
$\calD(L)$ is dense in $V^N$ \cite{K}. By abuse, we
do not distinguish in the notation $L$ from its restriction. 
We remark that $L^*$, the adjoint of $L$ is similarly
obtained from the coefficients
$\overline{a_{\beta\alpha}^{ji}}$.

Such an operator as holomorphic functional calculus. It
satisfies  resolvent estimates such as 
$$
\| (\zeta- ( L+\kappa))^{-1}\|_{op} \le dist(\zeta,
\Gamma)^{-1}, \zeta \notin \Gamma
$$
where $\Gamma$ is an open sector with vertex 0,  and half
angle
$\omega$ from the positive $x$-axis,
where
$w \in  [0, \pi/2)$ depends on $\lambda, \Lambda, N, m, n$.
Such estimates allow by Cauchy formula to define $f(L)$
for some appropriate holomorphic functions $f$ defined on
conic  neighborhoods of the spectrum of $L+\kappa$ (ie,
defined on larger open sectors).

In particular, one can take $f(\zeta)= \zeta^\alpha$ for
$\alpha\in [-1,1]$ and obtain the fractional powers of
$L+\kappa$. They are closed unbounded operators with the  
 expected properties such as 
$$
(L+\kappa)^{\alpha}(L+\kappa)^\beta=(L+\kappa)^{\alpha+\beta}
$$
when $\alpha+\beta \in [-1,1]$. In particular,
 $(L+\kappa)^{1/2}$  is the unique maximal-accretive square
root of
$L+\kappa$.

Kato first studied this question: is it possible to identify
the domains of the positive fractional powers of $L+\kappa$?

\section{In what square roots are critical?}

Kato found the following answer by abstract functional
analytic methods \cite{K1}. He proved that  for $\alpha \in
(0,1/2)$ then
$\calD((L+\kappa)^\alpha) =\calD(L^*+\kappa)^\alpha)$. This
result was completed by J.L. Lions \cite{Li} which found other
identifications by compex interpolation and one has
$\calD((L+\kappa)^\alpha) =
[L^2(\Omega), V^N]_{2\alpha}
$. Whenever such interpolation spaces are known then one gets
a result.

Also Lions proved that for any $\alpha \in (0,1)$,
$\calD((L+\kappa)^\alpha)=[L^2(\Omega), \calD(L)]_{\alpha}$
but this result is in practice useless as we do not know the
domain of $L$.  This implies nevertheless that whenever
 $\calD((L+\kappa)^{1/2}) $ and $\calD((L^*+\kappa)^{1/2})$
are both contained in $V^N$ then the three spaces are the
same. 

But the methods break down at $\alpha=1/2$ and  the
result cannot be true by purely abstracts reasonnings as we
see in the next section.  The remaining question is the
following.

\begin{conjecture} [Kato square root problem] Does
$\calD((L+\kappa)^{1/2})
$ coincide with the domain of the form $Q$?
\end{conjecture}

One case is easy. When $L$ is
self-adjoint then 
$$\|(L+\kappa)^{1/2} f\|^2_2 = \langle (L+\kappa)f, f\rangle
= Q(f,f) + \kappa\|f\|^2_2\ge\lambda\|\nabla^m f\|_2^2
$$
for all $f \in \calD(L)$. Thus,  $\calD((L+\kappa)^{1/2}) $ is
contained in $V^N$, hence the spaces coincide.

Let us see why $\alpha=1/2$ is critical. Let $n=1$, $m=1$ and
$N=1$. That is, consider  $L=DaD$ with $D={-i
d/dx}$ with domain $H^1(\RR)$ and $a$ is the multiplication by
a bounded {\bf real}-valued function
$a(x)$ on $\RR$ such that $ a \ge 1$. In such a case, $L$ is
self-adjoint and the domain of $L$ is the space of $f \in
H^1(\RR)$ such that $af' \in H^1(\RR)$. It is not too hard to
construct functions in the space [Actually, this space can
even be characterized by an adapted wavelet basis, see
\cite{AT1}].

By self-adjointness, we have
$\calD(L^{1/2})=H^1(\RR^n)$  [This holds for {\bf complex} $a$
with 
$\R a \ge 1$, but this is much harder].
Using interpolation we find that 
$$
\calD(L^{\alpha})= \begin{cases} H^{2\alpha}(\RR), &{\rm if}\
\alpha\in (0,1/2),\\
\{f\in H^1(\RR); af' \in H^{2\alpha-1}(\RR)\}, &{\rm if}\
\alpha\in (1/2,1). 
\end{cases}
$$ 

In one dimension, the onto-ness of ${-i d/dx}$ and the
1-1-ness of its adjoint make the understanding of the
domain of
$L$ easier. In higher dimensions, these properties are lost.

\section{Abstract methods are  unsufficient}

 We present an adaptation of an abstract counterexample by
McIntosh \cite{Mc2}. On
$H=
\ell^2(\ZZ)$, define an unbounded selfadjoint operator $D$ by
$De_j= 2^j e_j$ and a bounded operator
$B$ by $Be_j=\sum_{n
\in \ZZ} a_{n}e_{j+n}$, where $(e_j)$ is the
natural hilbertian basis of $H$ and $(b_n)$ is a
sequence of complex numbers such that $\hat
b(\theta)= \sum b_n e^{in\theta}$ satisfies
$\|\hat b\|_\infty=1$.  Clearly, the operator
$B$ has norm equal to $\|B\|=\|\hat b\|_\infty
=1$.
For $z \in \CC$ with $|z|<1$, one can define the
maximal-accretive operator
$L_z=DA_zD$ with   $A_z= Id
+ zB$  by the method of forms. Let $R_z= (L_z)^{1/2}$.

Assume that $\|R_z u\| \le c
\|Du\|$ for all $u \in \mathcal{D}(D)$ and
uniformly for $|z|
\le r <1$. As a function of
$z$,
$R_z$ is an operator valued holomorphic function
so that $R_0'D^{-1}$ is bounded on 
$H$. 
 Differentiating at $z=0$ the equation
$R_zR_z=L_z$, we find
$$
R_0'D+DR_0'=DBD.
$$
Solving for $R_0'$ one finds that 
$$
R_0'e_j=  2^j \sum c_n
e_{j+n}\ , \quad c_n= \frac{b_n 2^{n}} {1+2^n}.$$
Hence,  $\|R_0'D^{-1}\|= \|\hat c\|_\infty$
with evident notation. Now take $b_n = \frac{i}{
\pi n}$, then 
$\displaystyle \hat b(\theta)= - \frac{2}{\pi}
\sum_{n>0} \frac{\sin (n\theta) }{n}=\frac{\theta
}{
\pi} -1 $, $0 < \theta < 2\pi,$ so
that
$\|\hat b\|_\infty=1$. But
$\hat c(\theta)
\sim -\frac{i
}{
\pi}
\ln |\sin(\theta/2)\,|$ near 0 so that $\hat c$ is
not bounded. This is a contradiction, hence
$\|R_z u\|
\le c
\|Du\|$ fails for some $z$.

We shall find out that the Kato conjecture for elliptic
operators belongs to the realm of harmonic analysis.

\section{Why complex coefficients?}

Take two pure second order {\bf self-adjoint}  operators $L_1$
and
$L_2$ on $\RR^n$ defined as in Section  1  and denote by $A_1$
and
$A_2$ the matrix of coefficients corresponding to $L_1$ and
$L_2$. Is is true that 
\begin{equation} \label{regularity}
\|( L_1)^{1/2} f - ( L_2)^{1/2} f\|_2 \le C
\|A_1-A_2\|_\infty \|\nabla f\|_2\quad ?
\end{equation}
This apparently simple question is equivalent to 
asking about the  strong regularity of the (non-linear)
mapping
$${\rm coefficients}\
\mapsto {\rm square\ root}$$ from an open set in
$L^\infty(\Omega, E)$ into the space of bounded operators from
$V^N$ to
$L^2(\Omega,
\CC^N)$, where $E$ is some finite dimensional space. 

This question is highly non-trivial. The solution of the
conjecture for all possible {\bf complex} coefficients (or
least those  complex coefficients that are perturbations of
self-adjoint coefficients) gives us boundedness of this
mapping on complex balls, hence analyticity by the use of
complex function theory.

Here is an application of \eqref{regularity}. Consider the
solutions $u_k(t,x)$, $t\in \RR$, $x\in \RR^n$, $k=1,2$,
of the wave equations
$$\partial^2_t u_k(t) + L_ku_k(t)=0,\quad t\in \RR, x\in
\RR^n,
$$
with same Cauchy data $\partial_t u_{k_{\mid t=0}} = g \in
L^2(\RR^n)$ and $u_k(0)=f \in H^1(\RR^n)$. Then, starting
grom the ansatz
$$
u_k(t) =e^{it(L_k)^{1/2}}\tilde g + e^{-it(L_k)^{1/2}
}\tilde f
$$
and using \eqref{regularity} one obtains for $t>0$
$$
\|u_1(t)-u_2(t)\|_2 + \| \int_0^t \nabla (u_1(s)-u_2(s))\, ds
\|_2 \le C t\|A_1-A_2\|_\infty (\|\nabla f\|_2 + \|g\|_2).
$$
This estimate is sharp. It suffices to take $L_1=-\Delta$ and 
$L_2=-(1+b)\Delta$  with $b$ small to show this.

\section{The known results}

Here we state the positive answers to the above conjecture.

\begin{theorem} Let $n\ge 1$ and $L=-\dv A \nabla$ be a pure
second order
 operator on $\RR^n$. Then $\calD(L^{1/2})=H^1(\RR^n)$ with
the estimate $\|L^{1/2}f\|_2 \sim \|\nabla f\|_2$.
\end{theorem}

This is the result we shall
explain in the following sections.

The case $n=1$ was due to Coifman, McIntosh and Meyer in 1981
\cite{CMcM}. The general case is due to Hofmann, Lacey,
McIntosh, Tchamitchian and the author \cite{AHLMcT}. It came
after a succesful attempt in 2-d by Hofmann and McIntosh
(unpublished manuscript). See the introduction \cite{AHLMcT} 
for references to earlier partial results. 

\begin{theorem} Let $n\ge 1$ and $L$ be an homogeneous
elliptic  $N\times N$-system of arbitrary order $m$ on
$\RR^n$. Then $\calD((L+\kappa)^{1/2})=H^m(\RR^n,\CC^N)$ with
the estimate $\|(L+\kappa)^{1/2}f\|_2 \le C(\|\nabla^m
f\|_2^2+\kappa\|f\|_2^2)^{1/2} $. 
\end{theorem}

This result is due to Hofmann, McIntosh, Tchamitchian and the
author \cite{AHMcT}. 

\begin{theorem}\label{Omega} Let $n\ge 1$ and $L=-\dv A
\nabla$ be a pure second order
 operator on a proper open set $\Omega$ of $\RR$. Then one
has
$\calD(L^{1/2})=H^1(\Omega)$ with the estimate
$\|L^{1/2}f\|_2\le C(\|\nabla f\|_2+\|f\|_2)$ in the
following cases
\begin{itemize}
\item [(i)] $n=1$ and all possible choices of $\Omega$ and
$V$.
\item [(ii)] $n\ge 2$, $\Omega$ is a strongly Lipschitz
domain and $V=H^1_0(\Omega) $ (Dirichlet boundary condition)
or $V=H^1(\Omega)$ (Neumann boundary condition).
\end{itemize}
\end{theorem}

This theorem is due to Tchamitchian and the author.
In one dimension, this is achieved by  constructing an
adapted wavelet basis \cite{AT1}. We mention the approach by
interpolation methods and the result on $\RR$ by McIntosh,
Nahmod and the author \cite{AMcN1}. In higher dimensions, this
goes by transferring the result from $\RR^n$ \cite{AT4}. It
is likely that the method applies to second order systems
with Dirichlet or Neumann boundary conditions.

\begin{prop} Assume that $L$ is as in one of the
previous theorems. Then one can perturb $L$ by lower order
terms (ie, obtain an inhomogeneous operator) and still
answer positively the square root problem for the perturbed
operator. 
\end{prop}

We have separated this result from the others because it
is an ``abstract'' statement proved in
\cite {AT}, Chapter 0, Proposition 11. . Basically, any positive
result for the square root of a given homogeneous operator 
is ``stable'' under perturbations by lower order terms.

\section{Open problems} 

We list some problems ranked by level of difficulty, the
first being most likely more tractable. 

\begin{problem}
Find a  direct proof of Theorem \ref{Omega} 
following the ideas of 
\cite{AHLMcT}.
\end{problem}

\begin{problem} This problem was already posed by Lions.
Prove the Kato conjecture for second order operators with
mixed boundary conditions on strongly Lipschitz domains.
\end{problem}

\begin{problem} More generally,
prove the Kato conjecture for second order operators under
general boundary conditions on strongly Lipschitz domains.
\end{problem}

\begin{problem} Prove the Kato conjecture for higher order
operators or systems  with Dirichlet or
Neumann boundary conditions on smooth domains, then on
strongly Lipschitz domains. Study other types of boundary
conditions. 
\end{problem} 

\section{Harmonic analysis}

Our goal is to understand when a square function estimate
(SFE) of the form
\begin{eqnarray}
\label{eqn3:1}
\left(\int_0^\infty \|U_tf\|_2^2 \frac{dt }{
t}\right)^{1/2} \le C\|f\|_2, 
\end{eqnarray}
can hold,
where $(U_t)_{t>0}$ is a family of operators acting boundedly
and uniformly on $L^2(\RR^n)$.

We shall present the ideas in a model case and say how to
generalize them. Proofs will not be given and the reader is
referred \cite{AT} and \cite{CJ} for the $T(1)$ theorem. The
version of the $T(b)$ theorem given here is new. Related
ideas are in \cite{Ch}.

\subsection{The T(1) theorem}

The first part of the program is to find a simple statement
equivalent to SFE.

 We consider a model
case in which one can compute  $U_tf(x)$ as $$\int U_t(x,y)
f(y) dy$$  where the kernel
$U_t(x,y)$ is supported in
$|x-y|
\le t$ and satisfies
\begin{eqnarray}
\label{eqn3:13}
|U_t(x,y)| \le t^{-n}\qquad \hbox{ and }\qquad
|\nabla_y U_t(x,y)| \le t^{-n-1}.
\end{eqnarray}

Notice that only a regularity on the second variable is
imposed. 

Let $Q$ be a cube with side parallel to the axes. We denote
by 
 $|Q|$ its volume in $\RR^n$ and by $\ell(Q)$ its sidelength.
Also $cQ$ denotes the cube obtained by dilating $c$ times $Q$
from the centre of $Q$.
  If we apply \eqref{eqn3:1} to $f
={\bf 1}_{3Q}$ (the indicator function of $3Q$) and observe
 from \eqref{eqn3:13} that
$$
(U_t1)(x)=U_t({\bf 1}_{3Q})(x)$$
whenever  
$x\in Q$ and $0<t\le\ell(Q)$, then we obtain
$$\int_Q \! \int_0^{\ell(Q)}| (U_t1)(x)|^2 \,
\dxdtovert \le C|3Q|=C3^n|Q|.$$
Such an estimate means that $| (U_t1)(x)|^2 \,
\dxdtovert$ is a Carleson measure, that is a (positive
Borel regular) measure
$\mu$ on
$\RR^n \times (0,+\infty)$ such that
$$
\sup\frac{\mu(\calR_Q)}{|Q|} <+\infty
$$
where the supremum is taken over all cubes $Q$. We have set
$\calR_Q=Q\times (0,\ell(Q)]$. We denote this supremum by
$\|\mu\|_c$ and call it the Carleson norm of $\mu$.

There is a converse to this which begins with  the
celebrated Carleson inequality. 

\begin{lemma} Assume that $P_t$ is an
operator  with kernel similar to that of $U_t$ (only  a size
estimate are used at this point) then for any Carleson
measure
$\mu$, 
$$
\int_0^\infty\!\int_{\RR^n} |P_tf(x)|^2 d\mu(x,t) \le C
\|\mu\|_c \int_{\RR^n} |f|^2
$$
\end{lemma}

Assuming now that $| (U_t1)(x)|^2 \,
\dxdtovert$ is a Carleson measure, this tells us that 
the operator 
$$
f \mapsto (U_t1)\cdot (P_tf)
$$
satisfies  SFE. Hence,  the SFE for $U_t$ is the same as
the SFE for $V_t$ with
$$V_t=U_t -(U_t1)\cdot P_t.
$$ 
The latter operator has a kernel satisfying \eqref{eqn3:13}
(the regularity  for $P_t(x,y)$ in the second
variable is used here). If, in addition, we impose 
$$
P_t1=1
$$
 then we have  
 $$V_t1=0,$$ that is
\begin{equation}\label{cancellation}
\int_{\RR^n} V_t(x,y)\, dy=0.
\end{equation}
This cancellation condition permits  almost-orthogonality
arguments in a second step. 

Let us begin with the Schur Lemma.

\begin{lemma} Let $(\Delta_s)_{s>0}$   be a family of
self-adjoint (this is just to make life easy) operators on
$L^2(\RR^n)$ such that 
\begin{equation}\label{LP}
f= \int_0^\infty \Delta_s^2 f\, \frac{ds}{s}
\end{equation}
in the $L^2$-sense.
Assume also  the almost-orthogonality $L^2-L^2$ bound  
\begin{equation}\label{AO}
\|V_t\Delta_s||_{op}   \le C
\left(\inf\left(\frac{t}{s},
\frac{s}{t}\right)\right)^\alpha.
\end{equation}
for some $\alpha>0$.  Then $V_t$ satisfies SFE.
\end{lemma}

In practice, take $\Delta^*_t (=\Delta_t)$ with the similar
properties as the operator
$V_t$. Very often, $\Delta_t$ is an operator of convolution
type and \eqref{LP} is checked by use of the Fourier
transform. Now to see that the almost-orthogonality bound
holds we compute the kernel of $V_t\Delta_s$ as 
$$\int_{\RR^n} V_t(x,z) \Delta_s(z,y)\, dz.
$$
When $|x-y| \ge 2 \sup(t,s)$, then the support condition
gives us $0$, which is to say that the two functions of $z$
are orthogonal. When $|x-y| \le 2 \sup(t,s)$ then, we see
that  the function with smaller support oscillates while the
other  is regular on that support. Thus one can perform
an ``integration by parts'' by writing, if say $s\le t$,
$$
\int_{\RR^n} V_t(x,z) \Delta_s(z,y)\, dz=
\int_{\RR^n} (V_t(x,z)-V_t(x,y)) \Delta_s(z,y)\, dz.
$$
Using  the mean value inequality,  we get the bound 
$$
C\frac{s}{t} t^{-n}{\bf 1}_{|x-y|\le 2t}
$$
from which we obtain one of the almost-orthogonality bound.
The other one is exactly symmetric since we have the
cancellation condition \eqref{cancellation}.

Hence, the SFE for $V_t$ is always valid.
Let us summarize the results. 

\begin{theorem}[T(1) theorem] Assume $U_t$ and $P_t$ as above
with $P_t1=1$. Then, the followings are equivalent 
\begin{itemize}
\item [(i)] $U_t$ satisfies SFE. 
\item [(ii)]$(U_t1)\cdot P_t$ satisfies SFE. 
\item [(iii)]  $| (U_t1)(x)|^2 \,
\dxdtovert$ is a Carleson measure.
\end{itemize}
Moreover, one has
$$
\int_0^\infty\!\int_{\RR^n} |U_tf(x)-(U_t1)(x)\cdot
(P_tf)(x)|^2
\dxdtovert
\le C \int_{\RR^n} |f|^2.
$$
\end{theorem}

The idea of comparing $U_tf$ to $(U_t1)(P_tf)$ is natural in
probability where $U_t$ would be a positive linear operator.
It was brought into the topic of square function estimates
and Carleson measures by Coifman and Meyer \cite{CM1}.

\begin{remark} 

1) By handling tails, one can assume that $U_t(x,y)$ has some
integrable decay at infinity such as 
$$|U_t(x,y)| \le t^{\ep}(t+|x-y|)^{-n-\ep} , \quad \ep>0.$$
One can also replace the Lipschitz regularity by a H\"older
type regularity

2) One can take for $P_t$ a dyadic averaging operator: Given
a family of dyadic cube $Q$ of $\RR^n$, define 
$$
S_tf(x)= \frac{1}{|Q|} \int_Q f, \quad {\rm when}\ x\in Q\
{\rm and}\ 
\ell(Q)/2 < t
\le \ell(Q).$$
The difference is that the kernel of $S_t$ is not H\"older
smooth in its second variable. However, it is Sobolev
smooth, in the sense that it belongs to $H^{s}(\RR^n)$ when
$s\in (0,1/2)$. This is enough.

\end{remark} 

 In our applications, $U_t$ will neither have a nice
kernel, nor regularity in the second variable. 
Here is the statement which applies.

\begin{lemma}\label{lemmaLP} Let $U_t:L^2(\RR^n) \to L^2(\RR^n)$, $t>0$, be
a measurable family of bounded operators with $||U_t||_{op} \le 1$. Assume
that 
\begin{enumerate} 
\item[\rm (i)] $U_t$ has a kernel, $U_t(x,y)$, that is a measurable
function on $\RR^{2n}$ such that for some $m>n$ and for all $y\in \RR^n$
and $t>0$, 
$$
\int_{\RR^n}
\left(1+\frac{|x-y|}{t}\right)^{2m}|U_t(x,y)|^2\, dx \le
t^{-n}.
$$ 
\item[\rm (ii)] For any ball $B(y,t)$ with center at $y$ and radius $t$,
$U_t$ has a bounded extension from $L^\infty(\RR^n)$ to $L^2(B(y,t))$
with
$$
\frac{1}{t^n} \int_{B(y,t)} 
|U_tf(x)|^2\, dx \le \|f\|_\infty^2.
$$
 and  $U_t(f\calX_R)$
converges to $U_tf$ in $L^2(B(y,t))$  as $R\to\infty$  where 
$\calX_R$ stands for  the indicator function of the ball
$B(0,R)$.
\end{enumerate} 
Let $P_t$ be as above.    Then $U_tP_t$ 
satisfies SFE if and only if $| (U_t1)(x)|^2 \,
\dxdtovert$ is a Carleson measure.
Moreover, one has
$$
\int_0^\infty\!\int_{\RR^n} |U_tP_tf(x)-(U_t1)(x)\cdot
(P_tf)(x)|^2
\dxdtovert
\le C \int_{\RR^n} |f|^2.
$$
\end{lemma}

The idea of proof is to go back to the previous theorem by
using the operator $U_t^*U_tP_t$.

The same conclusion holds if one replaces $P_t$ by $S_t$.

\subsection{The T(b) theorem}

The next part of the program is to be able to obtain the
Carleson measure estimate involving $U_t1$. The ideas here
grew out from Semmes' work \cite{S}.

In practice, either $U_t1=0$ and there is nothing to do or
$U_t1\ne 0$ and it is usually impossible to compute. T(b)
theorems are useful tools designed to overcome such problems. 

Let us go back to a model operator $U_t$ as in the previous
section. Assume that for each cube $Q$, there are functions
$b_Q: 3Q \to \CC$ with the following properties
\begin{equation}\label{bq1}
\int_{3Q} |b_Q|^2 \le C|Q|,
\end{equation}
\begin{equation}\label{bq2}
|(S_tb_Q)(x)| \ge \delta, \quad {\rm for}\ (x,t) \in \calR_Q,
 \end{equation}
\begin{equation}\label{bq3}
(U_tb_Q)(x)=0 \quad {\rm for}\ (x,t) \in \calR_Q.
\end{equation}
The constant $C$ and $\delta$ are of course independent  of
$Q$. Here the dyadic cubes have been chosen so that $Q$ is
one of them. Then
\begin{align*}
\int_Q \! \int_0^{\ell(Q)}| (U_t1)(x)|^2 \,
\dxdtovert &\le \frac{1}{\delta^2} \int_Q \!
\int_0^{\ell(Q)}| (U_t1)(x)\cdot( S_tb_Q)(x)|^2 \,
\dxdtovert \\
 &\le \frac{2}{\delta^2} \int_Q \!
\int_0^{\ell(Q)}| (U_tb_Q)(x)|^2 \,
\dxdtovert \\
&\qquad + \frac{2}{\delta^2} \int_Q \!
\int_0^{\ell(Q)}| (V_tb_Q)(x)|^2 \,
\dxdtovert\\
&=\frac{2}{\delta^2} \int_Q \!
\int_0^{\ell(Q)}| (V_tb_Q)(x)|^2 \,
\dxdtovert\\
&\le C |Q|.
\end{align*}
The first inequality comes from \eqref{bq2}, the second from
the definition of $V_t$, then one uses \eqref{bq3} and the
last inequality comes from SFE for $V_t$ combined with 
\eqref{bq1}. 
 
Let us see how to relax the hypotheses. First,
\eqref{bq1} is OK as is. Secondly, \eqref{bq3} can clearly be
replaced by
\begin{equation}\label{bq3bis}
\int_Q \!
\int_0^{\ell(Q)}| (U_tb_Q)(x)|^2 \,
\dxdtovert
\le C |Q|.
\end{equation} 

Next, \eqref{bq2} implies in particular that $|b_Q(x)|\ge
\delta$ for $x \in Q$, which is often too strong. 
We shall need this lower bound only on a subset of
$\calR_Q$.

\begin{lemma} Let $\mu$ be a measure on $\RR^n\times
(0,\infty)$. Assume there are two constants $C>0$ and
$\eta\in (0,1)$ such that for each cube $Q$ one can find 
disjoint subcubes $Q_i$ of $Q$  with 
\begin{equation}\label{Qi}
\sum |Q_i| \le (1-\eta) |Q|
\end{equation}
and
$$\mu(\calR_Q\setminus \cup \calR_{Q_i}) \le C|Q|
$$
Then
$\mu$ is a Carleson measure and $\|\mu\|_c \le C/\eta$.
\end{lemma}

The proof is so simple that we give it.  Suppose a priori
that $\mu$ is a Carleson measure. We wish to obtain the
bound above.   Write
\begin{align*}\mu(\calR_Q) &=
\mu(\calR_Q\setminus \cup \calR_{Q_i}) + \sum
\mu( \calR_{Q_i})\\&
 \le
C|Q| + \|\mu\|_c \sum |Q_i|\\
&\le C|Q| + (1-\eta) \|\mu\|_c|Q|.
\end{align*}
It remains to divide by $|Q|$, to take the supremum over
$Q$ and to solve for $\|\mu\|_c$.

Thus one can replace \eqref{bq2} by 
\begin{equation}\label{bq2bis}
|S_tb_Q(x)| \ge \delta \quad {\rm for}\ (x,t) \in \calR_Q
\setminus \cup \calR_{Q_i}
\end{equation}
where the cubes $Q_i$ satisfy \eqref{Qi}. In the argument to
control $U_t1$, the LHS is only integrated on $\calR_Q
\setminus \cup \calR_{Q_i}$.
In other words, we allow a ``black hole'' region $\cup \calR_{
Q_i}$ on which we know nothing provided the ``bad'' cubes
$Q_i$ do not cover all of
$Q$.

Let me make a semantic digression. In French, a region
$\calR_Q$ is called ``fen\^etre de Carleson'', that is 
``Carleson window''. A very clean window lets the light
through. A window which may have some dark spots but not too
many of them still lets enough through. In other words, the
light goes through except for some ``black hole'' regions.

How to get the picture given by the ``lighted'' region
$\calR_Q
\setminus \cup \calR_{Q_i}$? The answer is by a
stopping-time argument.

The Carleson region $\calR_Q$ can be partitioned as the union
of  rectangles $$Q'\times ]\ell(Q')/2,\ell(Q')]$$ indexed by
all dyadic subcubes of $Q$ (they are called Whitney
rectangles), on which
$$(x,t)
\mapsto S_tb_Q(x)$$ is the constant function $$\frac{1}{Q'}
\int_{Q'} b_Q$$ (recall that $S_tb_Q(x)$ is  a dyadic average
of
$b_Q$ over a dyadic cube).   

Let us assume that $\int_Q b_Q=|Q|$. Let $\delta<1$. Consider
one of the dyadic children $Q'$ of $Q$, that is the cubes
obtained by subdividing
$Q$ with cubes with sidelength $\ell(Q)/2$.  We have two
options: 
\begin{itemize}
\item [(i)]  if the average gets too small, that is $$\R
\int_{Q'} b_Q
\le
\delta |Q'|,$$
then stop and select that cube.
\item [(ii)] otherwise
subdivide $Q'$ and argue similarly for each dyadic children.
\end{itemize}
 Keep going indefinitely
 and call $Q_i$
the cubes on which $b_Q$ has a small average. 

By construction, these cubes are disjoint, one can see
right away that the region $\calR_Q
\setminus \cup \calR_{Q_i}$ 
is the region where $\R (S_tb_Q)(x) \ge \delta$. 

It remains to see \eqref{Qi}.
Indeed, one has 
$$
\sum (1-\delta)|Q_i| \le \sum \R \int_{Q_i} 1-b_Q
= -\R\int_{Q\setminus \cup Q_i} 1-b_Q \le
C|Q|^{1/2}|Q\setminus \cup Q_i|^{1/2}
$$
by Cauchy-Schwarz inequality and \eqref{bq1}. One easily
concludes from there. Observe the crucial use of real parts 
in the above equality.

As we see, instead of asking for a pointwise lower bound 
$|b_Q|\ge \delta$ on $Q$, we only need a  lower bound on the 
average of $b_Q$ over $Q$, which is weaker.
 
Summarizing we have obtained the following theorem.

\begin{theorem}[local T(b) theorem] Let $U_t$ be as above.
Assume that one has a family of functions $b_Q:3Q\to \CC$ 
satisfying 
\eqref{bq1}, $|\int_Q b_Q| \ge |Q|$ and \eqref{bq3bis}, then
$| (U_t1)(x)|^2 \,
\dxdtovert$ is a Carleson measure. 
\end{theorem}

Again, one can state variations of the statement provided one
can make sense of $U_t1$ and have the SFE for $V_t$ or
$V_tS_t$.

\section{Back to square roots}

We are considering a pure second order operator 
$L=-\dv A \nabla$ with ellipticity constants $\lambda$ and
$\Lambda$ on $\RR^n$ ($\kappa=\kappa'=0$).

Since $L$ is maximal-accretive, a theorem of McIntosh and 
Yagi \cite{McY} asserts that 
$$
\|L^{1/2} f\|_2^2 \sim  \int_0^\infty \|(I+t^2L)^{-1}
tLf\|_2^2 \frac{dt}{t}.
$$
This can also be obtained using almost-orthogonality
arguments. If we set
$$
\theta_tF = (I+t^2L)^{-1} t\dv (AF)
$$
for $F=(F_1, \ldots, F_n)$ then we want to establish
\begin{equation}\label{SFEtheta}
\int_0^\infty \|\theta_t (\nabla f)\|_2^2 \frac{dt}{t} \le C
\|\nabla f\|_2^2.
\end{equation}
We are therefore facing a square function estimate and we
need to see what kind of estimates are available.

\subsection{Elliptic estimates} 

Pointwise bounds for the kernel of $\theta_t$ are 
false (Recall that we are merely assuming the coefficients of
$A$ to be measurable) even when the coefficients are real
(where the classical Aronson-De Giorgi-Nash-Moser theory can
be used). Moreover, this kernel will not be regular in its
second variable. 

In fact, there is no mathematical implication between the Kato
problem and pointwise bounds on heat kernels and vice-versa. 
The pointwise bounds are just handy when we have them.
 
What is possible to obtain are these
off-diagonal bounds in the mean.

\begin{lemma}\label{lemmaoffdiagonal} Let $E$ and $E_0$ be two
closed sets of
$\RR^n$ and set
$d=\dist(E,E_0)$, the distance between $E$ and $E_0$. Then
$$\int_E |(I+t^2L)^{-1}f(x)|^2\, dx  
\le C e^{-\tfrac{ d}{ct}} \int |f(x)|^2\, dx, \quad
\mathrm{Supp}\, f  \subset E_0, $$
$$\int_E |t\nabla(I+t^2L)^{-1}f(x)|^2\, dx  
\le C e^{-\tfrac{ d}{ct}} \int |f(x)|^2\, dx, \quad
\mathrm{Supp}\, f  \subset E_0, $$
$$\int_E |(I+t^2L)^{-1}t\dv (AF)\,(x)|^2\, dx  
\le C e^{-\tfrac{ d}{ct}} \int | F(x)|^2\, dx,
\quad
\mathrm{Supp}\, \, F\,  \subset E_0, $$
 where
 $c>0$ depends only on $\lambda$ and $\Lambda$, and $C$ on $n$, $\lambda$
and
$\Lambda$.
\end{lemma} 

These bounds will be sufficient for us thanks to the theory
developed for square function estimates. They are reminiscent
of the bounds found by Gaffney for Laplace-Beltrami operators
on manifolds.

These bounds also 
imply one can define in the $L^2_{loc}$ sense the resolvent
applied to functions with polynomial growth at infinity. In
particular, one has $$(I+t^2L)^{-1}(1)=1.$$ 

\begin{lemma}\label{lemmatechnicalestimates} For some $C$ depending only on
$n$, $\lambda$ and $\Lambda$, if $Q$ is a cube in $\RR^n$, 
$t\le
\ell(Q)$ and  $f$ is  Lipschitz function on $\RR^n$ then we have
$$
\int_{Q} |(I+t^2L)^{-1} f -f|^2 \le Ct^2 \|\nabla
f\|_\infty^2|Q|,
$$
$$
\int_{Q} |\nabla ((I+t^2L)^{-1} f -f)|^2 \le C\|\nabla
f\|_\infty^2|Q|.
$$ 
 \end{lemma}

\subsection{Applying the T(1) and T(b) theorems}

Choose $P_t$ to be here the  operator of convolution by 
$t^{-n}p(\frac{x}{t})$ with $\int p=1$ and $p\in
C^\infty_0(B(0,1))$, where $B(0,1)$ is the unit ball.

The first thing is to apply the theory of square functions 
in order to reduce to a Carleson measure estimate. 

We observe first that 
$$
(\theta_t - \theta_t P_t^2) (\nabla f)= (I -(I+t^2L)^{-1})
\frac{(I-P_t^2)f}{t} 
$$
so that  
$$
\int_0^\infty \|  (\theta_t - \theta_t P_t^2) (\nabla
f)\|_2^2
\frac{dt}{t} \le 4 \int_0^\infty \|  \frac{(I-P_t^2)f}{t}
\|_2^2
\frac{dt}{t} = 4 C\|\nabla f\|_2^2
$$
where the last equality follows from Plancherel's theorem. 

Now the elliptic estimates of Lemma
\ref{lemmaoffdiagonal} allows us to use Lemma \ref{lemmaLP} 
for $U_t=\theta_tP_t$.

Hence, SFE for $U_tP_t=\theta_tP_t^2$ is equivalent 
$|(\theta_t1)(x)|^2 \dxdtovert$ being a Carleson measure.
Here $1$ is the  $n\times n$ unit matrix.
Moreover, one can substitue $S_t$ for $P_t$.

Summarizing we see that \eqref{SFEtheta} reduces to the
proving that $|(\theta_t1)(x)|^2 \dxdtovert$ is a Carleson
measure. Moreover, one has 
\begin{equation}\label{CM}
\int_0^\infty\int_{\RR^n} |(\theta_t\nabla f)(x) -
(\theta_t1)(x) \cdot (S_t\nabla f)(x)|^2 \, \dxdtovert \le 
C\|\nabla f\|_2^2.
\end{equation}
Note that the product  $(\theta_t1)(x) \cdot (S_t\nabla
f)(x)$ is the dot product $u_1v_1 + \cdots + u_nv_n$ between
two vectors in $\CC^n$.

Now, we want to follow the ideas of the T(b) theorem.
There, the product was  over the complex field $\CC$. Since we
have now the dot product on $\CC^n$, we make a sectorial
decomposition of $\CC^n$. Let $\ep>0$ to be chosen later and cover $\CC^n$ with a
finite number depending on $\ep$ and $n$
 of cones $\calC_w$ associated to unit vectors  $w$   in
$\CC^n$ and defined by 
\begin{equation}
|u- (u | w)w| \le {\ep^{}}\ |(u| w)|.
\end{equation} 
Here $(\ |\ )$ is the complex inner product on $\CC^n$.
It suffices to   argue for each $w$ fixed  and to obtain a
Carleson measure estimate for 
$$\gamma_{t,w}(x)={\bf 1}_{\calC_w}((\theta_t1)(x))
(\theta_t1)(x),$$ where  
${\bf 1}_{\calC_w}$ denotes the indicator function of $\calC_w$. 

Fix $w$. We are looking for the analogs of the functions
$b_Q$. We call them $f_Q$. The requirements we are 
looking for 
are
\begin{equation} \label{fq1}
\int_{3Q} |\nabla f_Q|^2 \le C|Q|
\end{equation}
\begin{equation} \label{fq2}
|\int_Q \nabla f_Q| \ge \delta |Q|
\end{equation}
\begin{equation} \label{fq3}
\int_Q\!\int_0^{\ell(Q)} |(\theta_t\nabla f_Q)(x)|^2
\dxdtovert \le C |Q|
\end{equation}
and
\begin{equation} \label{fq4}
|\gamma_{t,w}(x)| \le C |\gamma_{t,w}(x) \cdot (S_t\nabla
f)(x) |
\end{equation} on ``good'' regions  $\calR_Q
\setminus \cup \calR_{Q_i}$ with not too many ``bad'' cubes
that is,
$\sum|Q_i|
\le (1-\eta)|Q|$.

The novelty is the last inequality which contains some
geometry. 

A 
candidate would be
$f_Q(x)= (x-x_Q|w)$ with $x_Q$ the centre of $Q$, because all
but the third inequality are fulfilled. Since
$\theta_t\nabla  = (I+t^2L)^{-1}tL $ it is natural to
approximate  $f_Q$ by applying  the resolvent to $f_Q$:
$$
f_Q^\ep= (I+\ep^2\ell(Q)^2L)^{-1} f_Q
$$
where $\ep$ is our small parameter. Note that $f_Q^\ep$ is
an approximation to $f_Q$ at the scale of $Q$.
It is defined on all of $\RR^n$ and 
Lemma  \ref{lemmatechnicalestimates} gives us
$L^2(3Q)$- estimates for  $f_Q-f_Q^\ep$ and its gradient. 

Hence, we obtain immediately
\eqref{fq1} and $C$ does not depend on $\ep$. We have
$$
\theta_t\nabla f_Q^\ep = (I+t^2L)^{-1}
\frac{t}{\ep^2\ell(Q)^2} (f_Q - f_Q^\ep)
$$
and  we deduce   
\eqref{fq3}. 

Now, to see \eqref{fq2} we observe that $\nabla f_Q=w^*$ (the
conjugate of $w$) and write 
$$
|\int_Q \nabla f_Q^\ep| \ge  \R (w^*|\int_Q \nabla f_Q^\ep)
=|Q| - \R \int_Q (w^*| \nabla f_Q- f_Q^\ep).  
$$
The inequality 
$$
\left|\int_Q \nabla h\,\right| \le 
{C}{\ell(Q)^{\frac{n-1}{2}}}\left(\int_Q  |h|^2\right)^{1/4} \left(\int_Q
|\nabla h|^2\right)^{1/4} 
$$
and Lemma 
\ref{lemmatechnicalestimates} imply
$$
\R \int_Q (w^*|\nabla f_Q- f_Q^\ep) \le C\ep^{1/2}|Q|
$$
and  \eqref{fq2} follows  provided $\ep$ is small enough.

It remains to obtain \eqref{fq4}. The stopping-time argument
of Section 7.2  would give us a lower bound of 
$\R (w^*|(S_t\nabla f_Q^\ep)(x))$ for $(x,t)$ in the
``good'' region. 
Given the fact that $\gamma_{t,w}(x)$ belongs to the cone
$\calC_w$ this is not enough. We also need to control
$|S_tf_Q^\ep(x)|$ on this ``good'' region. This means that
we have to introduce in the stopping-time argument a second
condition: starting from $Q$, we subdivide $Q$ dyadically
and stop the first time that either $\R \int_{Q'} (w^*|\nabla
f_Q) \le \delta|Q'|$ or $| \int_{Q'} \nabla  
f_Q | \ge C\ep^{-1}|Q'|$ where $C$ is appropriately chosen. 
As before, the union of the selected bad cubes cannot cover
all of
$Q$ if $\ep$ is small enough and we are done.  For details,
see
\cite{AHLMcT}.

\end{document}